\documentclass[11.05pt,a4paper]{amsart}
\usepackage{amssymb,amsmath}
\usepackage{latexsym}
\usepackage{amsthm,amsfonts,amssymb,mathrsfs}
\usepackage{rotating}
\usepackage[leqno]{amsmath}
\usepackage{xspace}
\usepackage[all]{xy}
\usepackage{longtable}

\headsep=1cm \numberwithin{equation}{section}
\newtheorem{theorem}{Theorem}[section]
\newtheorem{lemma}[theorem]{Lemma}

\newtheorem{corollary}[theorem]{Corollary}

\theoremstyle{definition}

\theoremstyle{remark}
\newtheorem{remark}[theorem]{Remark}

\newcommand{\Ass}{\operatorname{Ass}}

\newcommand{\Spec}{\operatorname{Spec}}

\newcommand{\Ht}{\operatorname{ht}}
\newcommand{\id}{\operatorname{id}}

\newcommand{\Gid}{\operatorname{Gid}}

\newcommand{\E}{\operatorname{E}}

\newcommand{\Ext}{\operatorname{Ext}}
\newcommand{\Supp}{\operatorname{Supp}}

\newcommand{\Tor}{\operatorname{Tor}}
\newcommand{\Hom}{\operatorname{Hom}}

\newcommand{\Ann}{\operatorname{Ann}}

\newcommand{\depth}{\operatorname{depth}}
\newcommand{\width}{\operatorname{width}}

\newcommand{\im}{\operatorname{im}}

\newcommand{\lo}{\longrightarrow}
\newcommand{\fm}{\frak{m}}
\newcommand{\fp}{\frak{p}}
\newcommand{\fq}{\frak{q}}

\newenvironment{prf}[1][Proof]{\begin{proof}[\bf #1]}{\end{proof}}
\begin{document}

\author[K. Divaani-Aazar, M. Nikkhah Babaei and M. Tousi]{Kamran Divaani-Aazar, Massoumeh Nikkhah Babaei
and Massoud Tousi}

\title[A criterion for dualizing modules]
{A criterion for dualizing modules}

\address{K. Divaani-Aazar, Department of Mathematics, Alzahra University, Vanak, Post Code 19834, Tehran,
Iran-and-School of Mathematics, Institute for Research in Fundamental Sciences (IPM), P.O. Box 19395-5746,
Tehran, Iran.}
\email{kdivaani@ipm.ir}

\address{M. Nikkhah Babaei, Department of Mathematics, Alzahra University, Vanak, Post Code 19834, Tehran,
Iran.}
\email{massnikkhah@yahoo.com}

\address{M. Tousi, Department of Mathematics, Shahid Beheshti University, G.C., Tehran, Iran-and-School of
Mathematics, Institute for Research in Fundamental Sciences (IPM), P.O. Box 19395-5746,
Tehran, Iran.} \email{mtousi@ipm.ir}


\keywords {$C$-injective modules; dualizing modules; semidualizing modules; trivial extensions.\\
The research of the first and third authors are supported by grants from IPM (No. 92130212 and No. 92130211,
respectively).}

\begin{abstract} We establish a characterization of dualizing modules among semidualizing modules. Let $R$
be a finite dimensional commutative Noetherian ring with identity and $C$ a semidualizing $R$-module.  We
show that $C$ is a dualizing $R$-module if and only if $\Tor_i^R(E,E')$ is $C$-injective for all $C$-injective
$R$-modules $E$ and $E'$ and all $i\geq 0.$
\end{abstract}

\maketitle

\section{Introduction}

Throughout this paper, $R$ will denote a commutative Noetherian ring with non-zero identity. The injective
envelope of an $R$-module $M$ is denoted by $\E_R(M)$.

A finitely generated $R$-module $C$ is called {\em semidualizing} if the homothety map $R\lo \Hom_R(C,C)$
is an isomorphism and $\Ext_R^i(C,C)=0$ for all $i>0$. Immediate examples of such modules are free
$R$-modules of rank one. A semidualizing $R$-module $C$ with finite injective dimension is called
{\em dualizing}. Although $R$ always possesses a semidualizing module, it does not possess a dualizing
module in general. Keeping \cite[Theorem 3.3.6]{BH} in mind, it is straightforward to see that the ring
$R$ possesses a dualizing module if and only if it is Cohen-Macaulay and it is homomorphic image of a finite 
dimensional Gorenstein ring.

Let $(R,\fm,k)$ be a local ring. There are several characterizations in the literature for a semidualizing
$R$-module $C$ to be dualizing. For instance, Christensen \cite[Proposition 8.4]{C} has shown that a
semidualizing $R$-module $C$ is dualizing if and only if the Gorenstein dimension of $k$ with respect to
$C$ is finite. Also, Takahashi et al. \cite[Theorem 1.3]{TYY} proved that a semidualizing $R$-module $C$
is dualizing if and only if every finitely generated $R$-module can be embedded in an $R$-module of finite 
$C$-dimension. Our aim in this paper is to give a new characterization for a semidualizing $R$-module $C$
to be dualizing.

Let $C$ be a semidualizing $R$-module. An $R$-module $M$ is said to be $C$-{\em projective} (respectively
$C$-{\em flat}) if it has the form $C\otimes_RU$ for some projective (respectively flat) $R$-module $U$.
Also, a $C$-{\em free} $R$-module is defined as a direct sum of copies of $C$. We can see that every
$C$-projective $R$-module is a direct summand of a $C$-free $R$-module and over a local ring every finitely
generated $C$-flat $R$-module is $C$-free. Also, an $R$-module $M$ is said to be $C$-{\em injective} if it
has the form $\Hom_R(C,I)$ for some injective $R$-module $I$.

Yoneda raised a question of whether the tensor product of injective modules is injective. Ishikawa in \cite
[Theorem 2.4]{I} showed that if $\E_R(R)$ is flat, then $E\otimes_RE'$ is injective for all injective
$R$-modules $E$ and $E'$. Further, Enochs and Jenda \cite[Theorem 4.1]{EJ} proved that $R$ is Gorenstein
if and only if for every injective $R$-modules $E$ and $E'$ and any $i\geq 0$, $\Tor_i^R(E,E')$ is injective.
We extend this result in terms of a semidualizing $R$-module. More precisely, for a semidualizing $R$-module
$C$, we show that the following are equivalent (see Theorem 2.7):

\begin{enumerate}
\item[(i)] $C_{\fp}$ is a dualizing $R_{\fp}$-module for all $\fp\in \Spec R$.
\item[(ii)] For any prime ideal $\fp$ of $R$ and any $i\geq 0$, $$\Tor_i^R(\E_C(R/\fp),\E_C(R/\fp))=
\begin{cases} 0&if \  \  i\neq \dim_{R_{\fp}}C_{\fp}\\
\E_C(R/\fp)&if \  \ i=\dim_{R_{\fp}}C_{\fp},\end{cases}$$
where $\E_C(R/\fp):=\Hom_R(C,\E_R(R/\fp))$.
\item[(iii)] For any $C$-injective $R$-modules $E$ and $E'$ and any $i\geq 0$, $\Tor_i^R(E,E')$
is $C$-injective.
\end{enumerate}

\section{The Results}

Let $\fp$ be a prime ideal of $R$. Recall that an $R$-module $M$ is said to have property $t(\fp)$ if for each
$r\in R-\fp$, the map $M\overset{r}\lo M$ is an isomorphism and if for each $x\in M$ we have $\fp^mx=0$ for some
$m\geq 1$. If an $R$-module $M$ has $t(\fp)$-property, then it has the structure as an $R_{\fp}$-module.
It is known that $\E_R(R/\fp)$ has $t(\fp)$-property.

To prove Theorem \ref{38}, which is our main result, we shall need the following five preliminary lemmas.

\begin{lemma}\label{31} Let $C$ be a semidualizing $R$-module. Then the following statements hold true.
\begin{enumerate}
\item[(i)] $\E_C(R/\fp):=\Hom_R(C,\E_R(R/\fp))$ has t$(\fp)$-property for each $\fp\in\Spec R$.
\item[(ii)] If $\fp$ and $\fq$ are two distinct prime ideals of $R$, then $\Tor_i^R(\E_C(R/\fp),\E_C(R/\fq))=0$
for all $i\geq 0$.
\end{enumerate}
\end{lemma}

\begin{prf} (i) As $\E_R(R/\fp)$ has t$(\fp)$-property, one can easily check that for any finitely generated
$R$-module $M$, the $R$-module $\Hom_R(M,\E_R(R/\fp))$ has $t(\fp)$-property.

(ii) By (i) $\E_C(R/\fp)$ has $t(\fp)$-property and  $\E_C(R/\fq)$ has $t(\fq)$-property. So, \cite[5]{EH}
implies that $$\Tor_i^R(\E_C(R/\fp),\E_C(R/\fq))=0$$ for all $i\geq 0$.
\end{prf}

\begin{lemma}\label{32} Let $(R,\fm,k)$ be a local ring, $C$ a semidualizing $R$-module and $I$ an Artinian
$C$-injective $R$-module. Then $\Hom_R(I,\E_R(k))$ is a finitely generated $\widehat{C}$-free $\widehat{R}
$-module.
\end{lemma}

\begin{prf} Denote the functor $\Hom_R(-,\E_R(k))$ by $(-)^{\vee}$. We have $I=\Hom_R(C,I')$ for some injective
$R$-module $I'$. Clearly, $C\otimes_RI$ is also an Artinian $R$-module. Since $$C\otimes_RI\cong C\otimes_R
\Hom_R(C,I')\cong \Hom_R(\Hom_R(C,C),I')\cong I',$$ we deduce that $I'$ is also Artinian. So, $I'\cong \overset{n}
\oplus \E_R(k)$ for some nonnegative integer $n$.

Now, one has$$I^{\vee}=\Hom_R(C,I')^{\vee}\cong C\otimes_RI'^{\vee}\cong \overset{n}\oplus \widehat{C},$$  and so
$I^{\vee}$ is a finitely generated  $\widehat{C}$-free $\widehat{R}$-module.
\end{prf}

In the next result, we collect some useful known properties of semidualizing modules. We may use them without
any further comments.

\begin{lemma}\label{34} Let $C$ be a semidualizing $R$-module and $\underline{r}:=r_1,\ldots, r_n$ a sequence of 
elements of $R$. The following statements hold.
\begin{enumerate}
\item[(i)] $\Supp_RC=\Spec R,$ and so $\dim_RC=\dim R$.
\item[(ii)] If $R$ is local, then $\widehat{C}$ is a semidualizing $\widehat{R}$-module.
\item[(iii)] $\underline{r}$ is a regular $R$-sequence if and only if $\underline{r}$ is a regular $C$-sequence.
\item[(iv)] If $\underline{r}$ is a regular $R$-sequence, then $C/( \underline{r} ) C$ is
a semidualizing $R/( \underline{r} )$-module.
\item[(v)] If $R$ is local and $\underline{r}$ is a regular $R$-sequence, then $C$ is a dualizing $R$-module if
and only if $C/( \underline{r} ) C$ is a dualizing $R/( \underline{r})$-module.
\end{enumerate}
\end{lemma}

\begin{prf} (i) and (ii) follow easily by the definition of a semidualizing module.

(iii) and (iv) are hold by \cite[Corollary 3.3.3]{S}.

(v) Assume that $R$ is local and $\underline{r}$ is a regular $R$-sequence. Then by (iv),
$C/(\underline{r}) C$ is a semidualizing $R/( \underline{r} )$-module.  On the other hand,
\cite[Corollary 3.1.15]{BH} yields that $$\id_{\frac{R}{(\underline{r})}}\frac{C}{(\underline{r})C}=
\id_RC-n.$$ This implies the conclusion.
\end{prf}

In the proof of the following result,  $R\ltimes C$ will denote the trivial extension of $R$ by $C$. For
any $R\ltimes C$-module $X$, its Gorenstein injective dimension will be denoted by $\Gid_{R\ltimes C}X$.
Also, we recall that for a module $M$ over a local ring $(R,\fm,k)$, the width of $M$ is defined by $\width_RM:
=\inf\{i\in \mathbb{N}_0|\Tor_i^R(k,M)\neq 0\}.$

\begin{lemma}\label{35} Let $(R,\fm,k)$ be a local ring and $C$ a semidualizing $R$-module. Then $\E_C(k)
\otimes_R\E_C(k)$ is a non-zero $C$-injective $R$-module if and only if $C$ is a dualizing $R$-module
of dimension $0$.
\end{lemma}

\begin{prf} Suppose that $\E_C(k)\otimes_R\E_C(k)$ is a non-zero $C$-injective $R$-module. As $\E_C(k)$
is Artinian, by \cite[Corollary 3.9]{KLS} the length of $\E_C(k)\otimes_R \E_C(k)$ is finite. So, also, 
$(\E_C(k)\otimes_R\E_C(k))^{\vee}$ has finite length. Since $$\Hom_R(\E_C(k),\widehat{C})\cong (\E_C(k)
\otimes_R \E_C(k)) ^{\vee},$$ by Lemma 2.2, we deduce that $\Hom_R(\E_C(k),\widehat{C})$ is isomorphic
to a direct sum of finitely many copies of $\widehat{C}$. This, in particular, implies that $\widehat{C}$
has finite length. Thus Lemma \ref{34} yields that $$\dim R =\dim_RC=\dim_{\widehat{R}}\widehat{C}=0,$$
and so, in particular, $R$ is complete. Next, one has
$$\begin{array}{lll}
\Hom_R(\E_C(k),R)&\cong \Hom_R(\E_C(k),\Hom_R(C,C))\\
&\cong \Hom_R(C,\Hom_R(\E_C(k),C))\\
&\cong \overset{n}\oplus \Hom_R(C,C)\\
&\cong R^n
\end{array}$$
for some $n>0$. This, in particular, implies that $$\Ann_R(\Hom_R(\E_C(k),R))=\Ann_RR.$$ Since $R$ is Artinian,
$\fm^t=0$ and $\fm^{t-1}\neq 0$ for some $t>0$. If for every $f\in \Hom_R(\E_C(k),R)$, $\im f\subseteq \fm$,
then $\fm^{t-1}f=0$ so $\fm^{t-1}\Hom_R(\E_C(k),R)=0$ a contradiction. Thus there is an epimorphism $\E_C(k)
\rightarrow R\rightarrow 0,$ and so $R$ is a direct summand of $\E_C(k)$. Next, \cite[Lemma 2.6]{HJ1} implies
that $R$ is a Gorenstein injective $R\ltimes C$-module. This yields that $C$ is a dualizing $R$-module, because
by \cite[Proposition 4.5]{HJ2}, one has $$\id_RC\leq \Gid_{R\ltimes C}R+\width_RR.$$

Conversely, if $C$ is a dualizing $R$-module of dimension 0, then $\dim R=0$ by Lemma \ref{34} (i). Hence,
$\E_R(k)$ is a dualizing $R$-module, and then by \cite[Theorem 3.3.4 (b)]{BH} we have $C\cong \E_R(k)$. Thus
$$\begin{array}{lll}
\E_C(k)\otimes_R \E_C(k)&\cong \Hom_R(\E_R(k),\E_R(k))\otimes_R\Hom_R(\E_R(k),\E_R(k))\\
&\cong R\otimes_R R\\
&\cong R\\
&\cong \Hom_R(C,\E_R(k)),
\end{array}$$
which is a non-zero $C$-injective $R$-module.
\end{prf}

\begin{remark}\label{36} (See \cite[(2.5)]{B}.) Let $M$ be an $R$-module and let $r\in R$ be a non-unit which
is  a non-zero divisor of both $R$ and $M$. Let $0 \rightarrow M \rightarrow I^0 \overset{d^0}\rightarrow I^1 
\rightarrow \cdots$ be a minimal injective resolution of $M$. Then there is a natural $R/(r)$-isomorphism
$M/(r)M\cong \Hom_R(R/(r),\im d^0)$ and $$0\rightarrow \Hom_R(R/(r),I^1)\rightarrow \Hom_R(R/(r),I^2)\rightarrow 
\cdots$$ is a minimal injective resolution of the $R/(r)$-module $M/(r)M$.
\end{remark}

Next, we recall the definition of the notion of co-regular sequences. Let $X$ be an $R$-module. An element $r$
of $R$ is said to be {\em  co-regular} on $X$ if the map $X\overset{r}\lo X$ is surjective. A sequence $r_1,
\dots , r_n$ of elements of $R$ is said to be a {\em  co-regular sequence} on $X$ if $r_i$ is co-regular on 
$(0:_M(r_1,\dots , r_{i-1}))$ for all $i=1,\dots , n$.

The following result plays a crucial role in the proof of Theorem \ref{38}.

\begin{lemma}\label{37} Let $(R,\fm,k)$ be a local ring and $C$ a semidualizing $R$-module. Let $r\in \fm$
be a non-zero divisor of $R$. Assume that $r$ is co-regular on $\Tor_i^R(\E_C(k),\E_C(k))$ for all $i.$ Then
for any $i\geq 0,$ we have a natural $\bar R$-isomorphism $$\Tor_{i-1}^{\bar R}(\E_{\bar C}(k),\E_{\bar C}(k))
\cong \Hom_R(\bar R,\Tor_i^R(\E_C(k),\E_C(k))),$$ where $\bar R:=R/( r ),$ $\bar C:=C/(
r) C,$ $\E_{C}(k):=\Hom_R(C,\E_{R}(k))$ and $\E_{\bar C}(k):=\Hom_{\bar R}(\bar C,\E_{\bar R}(k)).$
\end{lemma}

\begin{prf} Let $0\rightarrow I^0\rightarrow I^1\rightarrow\cdots $ be a minimal injective resolution of $C$.
Then $$\cdots\rightarrow\Hom_R(I^1,\E_R(k))\rightarrow\Hom_R(I^0,\E_R(k))\rightarrow 0$$ is a flat resolution
of $\E_C(k).$ Applying $\E_C(k)\otimes_R-$, we get the complex $$\cdots \rightarrow \E_C(k)\otimes_R\Hom_R
(I^1,\E_R(k))\rightarrow \E_C(k)\otimes_R\Hom_R(I^0,\E_R(k))\rightarrow 0.$$ We will denote $\E_C(k)\otimes_R
\Hom_R(I^i,\E_R(k))$ by $X_i$ and set $$X_{\bullet}:=\cdots \lo X_i\lo \cdots \lo X_1\lo X_0\lo 0.$$ Then for
each $i\geq 0,$ we have $H_i(X_{\bullet})=\Tor_i^R(\E_C(k),\E_C(k))$.

By Remark \ref{36}, $$0\rightarrow \Hom_R(\bar R,I^1)\rightarrow \Hom_R(\bar R,I^2)\rightarrow\cdots$$ is a
minimal injective resolution of $\bar C$ as an $\bar R$-module. So, $$\cdots \rightarrow\Hom_{\bar R}(\Hom_R(
\bar R,I^2),\E_{\bar R}(k))\rightarrow\Hom_{\bar R}(\Hom_R(\bar R,I^1),\E_{\bar R}(k))\rightarrow 0$$ is a flat 
resolution of
$\E_{\bar C}(k)$ as an $\bar R$-module. Thus for each $i\geq 1,$ the $\bar R$-module $\Tor_{i-1}^{\bar R}
(\E_{\bar C}(k),\E_{\bar C}(k))$ is isomorphic to the $i$th homology of the following complex $$(\star)
\cdots \lo \E_{\bar C}(k)\otimes_{\bar R}\Hom_{\bar R}(\Hom_R(\bar R,I^2),\E_{\bar R}(k))\lo \E_{\bar C}(k)
\otimes_{\bar R}\Hom_{\bar R}(\Hom_R(\bar R,I^1),\E_{\bar R}(k))\rightarrow 0.$$ We shall show that
the later complex is isomorphic to the complex $Y_{\bullet}:=\Hom_R(\bar R,X_{\bullet}).$

Noting that $\E_{\bar R}(k)\cong \Hom_R(\bar R,\E_R(k))$ and using Adjointness yields that $$\E_{\bar C}(k)
=\Hom_{\bar R}(\bar C,\E_{\bar R}(k))\cong \Hom_R(\bar R,\E_C(k)).$$ Hence for each $i\geq 0$, by using
Adjointness, Hom-evaluation and Tensor-evaluation, one has the
following natural $\bar R$-isomorphisms:
$$\begin{array}{lll}
\E_{\bar C}(k)\otimes_{\bar R}\Hom_{\bar R}(\Hom_R(\bar R,I^i),\E_{\bar R}(k))&\cong
\E_{\bar C}(k)\otimes_{\bar R}\Hom_{\bar R}(\Hom_R(\bar R,I^i),\Hom_R(\bar R,\E_R(k)))\\
&\cong \E_{\bar C}(k)\otimes_{\bar R}\Hom_R(\Hom_R(\bar R,I^i),\E_R(k))\\
&\cong \E_{\bar C}(k)\otimes_{\bar R}(\bar R\otimes_R\Hom_R(I^i,\E_R(k)))\\
&\cong \Hom_R(\bar R,\E_C(k))\otimes_R\Hom_R(I^i,\E_R(k))\\
&\cong \Hom_R(\bar R,\E_C(k)\otimes_R\Hom_R(I^i,\E_R(k)))\\
&\cong Y_i.
\end{array}$$
Note that $\Hom_R(I^i,\E_R(k))$ is a flat $R$-module. As $r$ is  a non-zero divisor of $R$, it is also a
non-zero divisor of $C$. This implies that  $r$ is a non-zero divisor of $I^0,$ and so $\Hom_R(\bar R,I^0)=0.$
Thus  $$Y_0\cong \E_{\bar C}(k)\otimes_{\bar R}\Hom_{\bar R}(\Hom_R(\bar R,I^0),\E_{\bar R}(k))=0.$$ Therefore,
the two complexes $(\star)$ and $Y_{\bullet}$ are isomorphic, and so we deduce that $\Tor_{i-1}^{\bar R}
(\E_{\bar C}(k),\E_{\bar C}(k))=H_i(Y_{\bullet})$ for all $i\geq 0$.

Since $r$ is a non-zero divisor of $C$, it is co-regular on $\E_C(k)$, and so it is co-regular on $X_i$
for all $i.$ Thus, we can deduce the following exact sequence of complexes $$0\lo Y_{\bullet}\lo X_{\bullet}
\overset{r}\lo X_{\bullet} \lo 0.$$ It yields the following exact sequences of modules $$\cdots \lo 
\Tor_{i+1}^R(\E_C(k),\E_C(k))\overset{r}\lo \Tor_{i+1}^R(\E_C(k),\E_C(k))\lo \Tor_{i-1}^{\bar R}
(\E_{\bar C}(k),\E_{\bar C}(k))\overset{f_i}\lo \Tor_i^R(\E_C(k),\E_C(k))$$ $$\overset{r}\lo \Tor_i^R(\E_C(k),\E_C(k))
\lo \cdots.$$ As $r$ is a co-regular element on $\Tor_i^R(\E_C(k),\E_C(k))$ for all $i$, we deduce that
$f_i$ is a monomorphism for all $i$. This implies our desired isomorphisms.
\end{prf}

\begin{theorem}\label{38} Let $C$ be a semidualizing $R$-module. The following are equivalent:
\begin{enumerate}
\item[(i)]  $C_{\fp}$ is a dualizing $R_{\fp}$-module for all $\fp\in \Spec R$.
\item[(ii)] For any prime ideal $\fp$ of $R$ and any $i\geq0$,
$$\Tor_i^R(\E_C(R/\fp),\E_C(R/\fp))=
\begin{cases} 0&if \  \ i\neq \dim_{R_{\fp}}C_{\fp}\\
\E_C(R/\fp)&if \   \ i=\dim_{R_{\fp}}C_{\fp},\\
\end{cases}$$
where $\E_C(R/\fp):=\Hom_R(C,\E_R(R/\fp))$.
\item[(iii)] For any $C$-injective $R$-modules $E$ and $E'$ and any $i\geq0$, $\Tor_i^R(E,E')$
is $C$-injective.
\end{enumerate}
\end{theorem}

\begin{prf} (i)$\Rightarrow$(ii) Let $\fp$ be a prime ideal of $R$. There are natural $R_{\fp}$-isomorphisms 
$\E_C(R/\fp)\cong \E_{C_{\fp}}(R_{\fp}/\fp R_{\fp})$ and $$\Tor_i^R(\E_C(R/\fp),\E_C(R/\fp))\cong 
\Tor_i^{R_{\fp}}(\E_{C_{\fp}}(R_{\fp}/\fp R_{\fp}),\E_{C_{\fp}}(R_{\fp}/\fp R_{\fp}))$$ for all $i\geq 0.$
Hence, we can complete the proof of this part by showing that if $C$ is a dualizing module of a local ring
$(R,\fm,k)$, then $$\Tor_i^R(\E_C(k),\E_C(k))=\begin{cases} 0&i\neq\dim_RC\\ \E_C(k)&i=\dim_RC. \end{cases}$$

Set $d:=\dim_RC$. As $C$ is a dualizing $R$-module, \cite[Theorem 3.3.10]{BH} implies that for any prime ideal
$\fp,$ one has $$\mu^i(\fp,C)=\begin{cases}0&i\neq \Ht \fp\\ 1&i=\Ht \fp.\end{cases}$$ So, if $I^{\bullet}=0
\rightarrow I^0\rightarrow I^1\rightarrow \cdots$ is a minimal injective resolution of $C$, then $I^d\cong
\E_R(k)$ and for any $i\neq d$, $\E_R(k)$ is not a direct summand of $I^i$. In particular, $\Hom_R(R/\fm,I^i)=0$
for all $i\neq d$. Now, $\Hom_R(I^{\bullet},\E_R(k))$ is a flat resolution of $\E_C(k).$  Clearly, one has 
$$\E_C(k)\otimes_R \Hom_R(I^d,\E_R(k))\cong \E_C(k)\otimes_R\widehat{R}\cong \E_C(k).$$ Next, let $i\neq d$.
Since $\Hom_R(I^i,\E_R(k))$ is a flat $R$-module, \cite[Theorem 23.2 (ii)]{M} implies that $$\Ass_R(\E_C(k)
\otimes_R\Hom_R(I^i,\E_R(k)))=\Ass_R(R/\fm\otimes_R\Hom_R(I^i,\E_R(k))).$$
But, $$R/\fm\otimes_R\Hom_R(I^i,\E_R(k))\cong \Hom_R(\Hom_R(R/\fm,I^i),\E_R(k))=0,$$ and so $\E_C(k)
\otimes_R\Hom_R(I^i,\E_R(k))=0.$ Therefore, it follows that the complex $\E_C(k)\otimes_R\Hom_R(I^{\bullet},
\E_R(k))$ has $\E_C(k)$ in its $d$-place and $0$ in its other places. Thus, we deduce that $$\Tor_i^R(\E_C(k),
\E_C(k))=H_i(\E_C(k)\otimes_R \Hom_R(I^{\bullet},\E(k)))=
\begin{cases}
0&i\neq d\\\E_C(k)&i=d.
\end{cases}$$

(ii)$\Rightarrow$(iii) Let $E$ be an injective $R$-module. Since $E\cong \underset{\fp\in \Spec R} \bigoplus
\E_R(R/\fp)^{\mu^0(\fp,E)}$ and $C$ is finitely generated, we have $$\Hom_R(C,E)\cong \underset{\fp\in \Spec R} 
\bigoplus \E_C(R/\fp)^{\mu^0(\fp,E)}.$$ As $R$ is Noetherian, clearly any direct sum of $C$-injective $R$-modules
is again $C$-injective, and so (ii) yields (iii) by Lemma \ref{31} (ii).

(iii)$\Rightarrow$(i) It is easy to check that a given $R_{\fp}$-module $M$ is $C_{\fp}$-injective if and only
if it is the localization at $\fp$ of a $C$-injective $R$-module. Thus, it is enough to show that if $C$ is a
semidualizing  module of a local ring $(R,\fm,k)$ such that $\Tor_i^R(E,E')$ is $C$-injective for all $C$-injective
$R$-modules $E$ and $E'$ and all $i\geq 0,$ then $C$ is dualizing.

Let $\underline{r}=r_1,\ldots, r_d\in \fm$ be a maximal regular $R$-sequence. Then $\underline{r}$ is also a
regular $C$-sequence. It is easy to verify that $\underline{r}$ is a co-regular sequence on any $C$-injective 
$R$-module, and consequently $\underline{r}$ is a co-regular sequence on  $\Tor_i^R(\E_C(k),\E_C(k))$ for all
$i\geq 0$.  Letting $\bar R:=R/(\underline{r})$ and $\bar C:=C/(\underline{r}) C,$  by Lemma \ref{34} (iv),
it turns out that $\bar C$ is a semidualizing $\bar R$-module. Making repeated use of Lemma \ref{37}, we can
establish the following natural $\bar R$-isomorphism $$\E_{\bar C}(k)\otimes_{\bar R}\E_{\bar C}(k)\cong
\Hom_R(\bar R,\Tor_d^R(\E_C(k),\E_C(k))).$$ So, $\E_{\bar C}(k)\otimes_{\bar R}\E_{\bar C}(k)$ is a
$\bar C$-injective $\bar R$-module. Lemma \ref{34} implies that $$\depth_{\widehat{\bar R}}\widehat{\bar C}
=\depth_{\bar R}\bar C=\depth_{\bar R}\bar R=0,$$ and so there are natural inclusion maps $k\overset{i}
\hookrightarrow \bar C$ and $k\overset{j}\hookrightarrow \widehat{\bar C}$. By applying the functor
$\Hom_{\bar R}(-,\E_{\bar R}(k))$ on $i$, we get an epimorphism $\E_{\bar C}(k)\twoheadrightarrow k$.
Next, by applying the functor $\Hom_{\bar R}(-,\widehat{\bar C})$ on the later map, we see that
$$\Hom_{\bar R}(\E_{\bar C}(k)\otimes_{\bar R}\E_{\bar C}(k),\E_{\bar R}(k))\cong \Hom_{\bar R}(\E_{\bar C}(k),
\widehat{\bar C})\neq 0.$$
Hence, $\E_{\bar C}(k)\otimes_{\bar R}\E_{\bar C}(k)$ is a non-zero $\bar C$-injective $\bar R$-module, and
so Lemma \ref{35} yields that $\bar C$ is a dualizing $\bar R$-module. Now, by Lemma \ref{34} (v), we deduce
that $C$ is a dualizing $R$-module.
\end{prf}

We end the paper with the following immediate corollary.

\begin{corollary}\label{39} Let $R$ be a finite dimensional ring and $C$ a semidualizing $R$-module. Then
$C$ is a dualizing $R$-module if and only if $\Tor_i^R(E,E')$ is $C$-injective for all $C$-injective
$R$-modules $E$ and $E'$ and all $i\geq 0$.
\end{corollary}


\end{document}